\documentclass[11pt,leqno]{article} 
\usepackage{graphics}
\newtheorem{thm}{Theorem}[section]
\newtheorem{lma}{Lemma}[section]
\newtheorem{cor}{Corollary}

\newcommand{\beqa}{\begin{eqnarray}}
\newcommand{\eeqa}{\end{eqnarray}}

\newcommand{\pf}{\noindent {\bf Proof:} $\s$ }
\newcommand{\epf}{ \hfill$\diamondsuit$ \medskip}

\newcommand{\ds}{\displaystyle}
\newcommand{\beq}{\begin{equation}}
\newcommand{\eeq}{\end{equation}}
\newcommand{\lbl}{\label}
\newcommand{\s}{\; \;}

\newcommand{\la}{\lambda}
\newcommand{\mb}{\mbox}
\newcommand{\ra}{\rightarrow}
\newcommand{\al}{\alpha}
\title{Families of solution curves for some non-autonomous problems}

\author{
Philip Korman   \\ 
Department of Mathematical Sciences \\ 
University of Cincinnati \\ 
Cincinnati Ohio 45221-0025 \\
kormanp@ucmail.uc.edu
}

\date{}

\begin{document}

\maketitle
\begin{abstract} 
The paper  studies families of positive solution curves for non-autonomous two-point problems 
\[
u''+\la f(u)-\mu g(x)=0, \s -1<x<1, \s u(-1)=u(1)=0 \,,
\]
depending on two  positive parameters $\la$ and $\mu$. We regard $\la$ as a primary parameter, giving us the solution curves, while the secondary parameter $\mu$ allows for evolution of these curves.  We give conditions under which the solution curves do not intersect, and the maximum value of solutions provides a global parameter. Our primary application is to constant yield harvesting for diffusive logistic equation. We implement numerical computations of the solution curves, using continuation in a global parameter,
a technique that we developed in \cite{K1}. 
 \end{abstract}

\begin{flushleft}
\end{flushleft}

\begin{flushleft}
\end{flushleft}

\section{Introduction }
\setcounter{equation}{0}

We study positive solutions  of non-autonomous two-point problems 
\beq
\lbl{1}
u''+\la f(u)-\mu g(x)=0, \s -1<x<1, \s u(-1)=u(1)=0 \,,
\eeq
depending on two  positive parameters $\la$ and $\mu$. We assume that  $f(u) \in C^2( \bar R_+)$, and $g(x) \in  C^1(-1,1) \cap C[-1,1]$ satisfies
\beq
\lbl{2}
g(-x)=g(x), \s \mbox{for $x \in (0,1)$} \,,
\eeq
\beq
\lbl{3}
g(0)>0, \s \mbox{and} \s xg'(x) \geq 0, \s \mbox{for $x \in (-1,1)$} \,.
\eeq
In case $g(x)$ is a constant, one can use the time map method, see K.C. Hung and S.H. Wang \cite{W1}, \cite{W2} who have studied similar multiparameter problems, or the book by S.P. Hastings and J.B.  McLeod \cite{H}. We show that under the conditions (\ref{2}) and (\ref{3}) one can still get detailed results on the solution curves $u=u(x,\la)$, where we regard $\la$ as a primary parameter, and on the evolution of these curves when the secondary parameter $\mu$ changes. We say that  the solution curves $u=u(x,\la)$ are the {\em $\la$-curves}. We also consider the $\mu$-curves, by regarding $\la$ as the secondary parameter.
\medskip

By B. Gidas, W.-M. Ni and L. Nirenberg \cite{GNN}, any positive solution of (\ref{1}) is an even function, and moreover $u'(x)<0$ for $x \in (0,1)$. It follows that $u(0)$ is the maximum value of the solution $u(x)$. Our first result says that $u(0)$ is a {\em global parameter}, i.e., its value uniquely determines the solution pair $(\la,u(x))$ ($\mu$ is assumed to be fixed). It follows that a planar curve $(\la,u(0))$ gives a faithful representation of the solution set of (\ref{1}), so that  $(\la,u(0))$ describes the {\em global solution curve}. Then we show positivity of any non-trivial solution of the linearized problem for (\ref{1}). This allows us to compute the direction of turn for convex and concave $f(u)$.
\medskip

Turning to the secondary parameter $\mu$, we show that solution curves at different $\mu$'s do not intersect, which allows us to discuss the evolution of solution curves in $\mu$.  
\medskip

We apply our results to a logistic model with fishing. S.  Oruganti, J. Shi, and R. Shivaji \cite{S} considered a class of general elliptic equations on an arbitrary domain, which includes 
\[
u''+\la u(1-u)-\mu g(x)=0, \s -1<x<1, \s u(-1)=u(1)=0 \,.
\]
They proved that for $\la$ sufficiently close to the principal eigenvalue $\la _1$, the $\mu$-curves are as in the Figure $2$ below. We show that such curves are rather special, with the solution curves as in the Figure $3$ below being more common.
Our approach is  to study the $\la$-curves first, leading to the understanding of the $\mu$-curves. We obtain an exhaustive  result in case $g(x)$ is a constant. The parameter $\mu>0$ quantifies the amount of fishing in the logistic model. We also consider the case $\mu <0$, corresponding to ``stocking" of fish.
\medskip

Using the fact that $u(0)$ is a  global parameter, we implement numerical computations of the solution curves, illustrating our results.
We use continuation in a global parameter,
a technique that we developed in \cite{K1}.

\section{Families of solution curves }
\setcounter{equation}{0}
\setcounter{thm}{0}
\setcounter{lma}{0}

The following result is included in  B. Gidas, W.-M. Ni and L. Nirenberg \cite{GNN}, see also P. Korman \cite{K14} for an elementary proof.
\begin{lma}\lbl{lma:1}
Under the conditions (\ref{2}) and (\ref{3}), any positive solution of
 (\ref{1}) is an even function,
with $u'(x)<0$ for all $x \in (0,1]$, so that $x=0$ is a point of global maximum.
\end{lma}

We begin by considering the secondary parameter $\mu$ to be fixed. To stress that, we call $h(x)=\mu g(x)$, and consider positive solutions of 
\beq
\lbl{4}
u''+\la f(u)-h(x)=0, \s -1<x<1, \s u(-1)=u(1)=0 \,.
\eeq
\begin{lma}\lbl{lma:2}
Assume that $f(u) \in C(\bar R_+)$ satisfies $f(u)>0$ for $u>0$, and $h(x)$ satisfies the conditions (\ref{2}) and (\ref{3}).
Then $u(0)$, the maximum value of any positive solution, uniquely identifies the solution pair $(\la, u(x))$.
\end{lma}

\pf
Observe from (\ref{4}) that $f(u(0))>0$ for any  positive solution $u(x)$. Let $(\la _1, v(x))$ be another solution of (\ref{4}), with $v(0)=u(0)$, $v'(0)=u'(0)=0$, and $\la _1>\la$.  From the equation (\ref{4}), $v''(0)<u''(0)$, and hence $v(x)<u(x)$ for small $x>0$. Let $\xi \leq 1$ be their first point of intersection, i.e., $u(\xi)=v(\xi)$. Clearly
\beq
\lbl{5}
|v'(\xi)| \leq |u'(\xi)| \,.
\eeq
Multiplying the equation (\ref{4}) by $u'$, and integrating over $(0,\xi)$, we get
\[
\frac12 {u'}^2(\xi)+\la \left[ F(u(\xi))-F(u(0)) \right]-\int_0^{\xi} h(x) u'(x) \, dx=0 \,,
\]
where $F(u)=\int_0^u f(t) \, dt$. Integrating by parts, we conclude
\[
\frac12 {u'}^2(\xi)=\la \left[ F(u(0))-F(u(\xi)) \right]+h(\xi)u(\xi)-h(0)u(0)-\int_0^{\xi} h'(x) u(x) \, dx \,.
\]
Similarly,
\[
\frac12 {v'}^2(\xi)=\la _1 \left[ F(u(0))-F(u(\xi)) \right]+h(\xi)u(\xi)-h(0)u(0)-\int_0^{\xi} h'(x) v(x) \, dx \,,
\]
and then, subtracting,
\[
\frac12 \left[{u'}^2(\xi)-{v'}^2(\xi) \right]=(\la -\la _1) \left[ F(u(0))-F(u(\xi)) \right]+\int_0^{\xi} h'(x) (v(x)-u(x)) \, dx \,.
\]
Since $v(x)<u(x)$ on $(0,\xi)$, the second term on the right is non-positive, while the first term on the right is negative, since $F(u)$ is an increasing function. It follows that $ |u'(\xi)|<|v'(\xi)|$, which contradicts (\ref{5}).
\epf

\begin{lma}\lbl{lma:3}
Assume that $f(u) \in C(\bar R_+)$, and $h(x)$ satisfies the conditions (\ref{2}) and (\ref{3}).
Then the  curves of positive solutions  of (\ref{1}) in $(\la, u(0))$ plane, computed at different $\mu$'s, do not intersect.
\end{lma}

\pf
Assume, on the contrary, that $v(x)$ is a solution of 
\beq
\lbl{7}
v''+\la f(v)-\mu _1 g(x)=0, \s -1<x<1, \s v(-1)=v(1)=0 \,,
\eeq
with $\mu _1 >\mu$, but $u(0)=v(0)$, where $u(x)$ is a solution of (\ref{1}). Then $u''(0)<v''(0)$, and hence $u(x)<v(x)$ for small $x>0$. Let $\xi \leq 1$ be their first point of intersection, i.e., $u(\xi)=v(\xi)$. Clearly
\beq
\lbl{8}
|u'(\xi)| \leq |v'(\xi)| \,.
\eeq
Multiplying the equation (\ref{1}) by $u'$, and integrating over $(0,\xi)$, we get
\[
\frac12 {u'}^2(\xi)+\la \left[ F(u(\xi))-F(u(0)) \right]=\mu \int_0^{\xi} g(x) u'(x) \, dx \,.
\]
Similarly, using (\ref{7}), we get
\[
\frac12 {v'}^2(\xi)+\la \left[ F(u(\xi))-F(u(0)) \right]=\mu _1 \int_0^{\xi} g(x) v'(x) \, dx \,.
\]
Subtracting, we obtain
\beqa \nonumber
& \frac12 \left[{u'}^2(\xi)-{v'}^2(\xi) \right]=\mu \int_0^{\xi} g(x) u'(x) \, dx-\mu _1 \int_0^{\xi} g(x) v'(x) \, dx \\ \nonumber
& > \mu _1 \left[\int_0^{\xi} g(x) u'(x) \, dx-\int_0^{\xi} g(x) v'(x) \, dx \right] \\ \nonumber
& = \mu _1 \int_0^{\xi} g'(x) \left(v(x)-u(x) \right) \, dx>0. \nonumber
\eeqa
Hence, $ |u'(\xi)|>|v'(\xi)|$, which contradicts (\ref{8}).
\epf

\begin{cor}\lbl{cor:1}
Assume that $\la$ is fixed in (\ref{1}), and $\mu$ is the primary parameter.  Assume that $f(u) \in C(\bar R_+)$ satisfies $f(u)>0$ for $u>0$, and $g(x)$ satisfies the conditions (\ref{2}) and (\ref{3}).
Then the maximum value of solution $u(0)$ is a global parameter, i.e., it   uniquely identifies the solution pair $(\mu, u(x))$.
\end{cor}

\pf
If at some $\la _0$ we had another solution pair $(\mu _1,u_1(x))$ with $u(0)=u_1(0)$, then the $\la$-curves at $\mu$ and $\mu _1$ would intersect at $(\la _0,u(0))$, contradicting Lemma \ref{lma:2}.
\epf

The linearized problem for (\ref{1}) is 
\beq
\lbl{10}
w''+\la f'(u)w=0, \s -1<x<1, \s w(-1)=w(1)=0 \,.
\eeq
We call the solution of (\ref{1}) {\em singular} if (\ref{10}) has non-trivial solutions. Since the solution set of (\ref{10}) is one-dimensional (parameterized by $w'(-1)$),
it follows that $w(-x)=w(x)$, and $w'(0)=0$.

\begin{lma}\lbl{lma:4}
Assume that $f(u) \in C^1(\bar R_+)$, and $g(x)$ satisfies the conditions (\ref{2}) and (\ref{3}), and let $u(x)$ be a positive solution of (\ref{1}).
Then any non-trivial solution of (\ref{10}) is of one sign, i.e., we may assume that $w(x)>0$ for all $x \in (-1,1)$.
\end{lma}

\pf
Assuming the contrary, we can find a point $\xi \in (0,1)$ such that $w(\xi)=w(1)=0$, and $w(x)>0$ on $(\xi,1)$. Differentiate  the equation (\ref{1})
\[
u'''+\la f'(u)u'-\mu g'(x)=0 \,.
\]
Combining this with (\ref{10}),
\beq
\lbl{11}
\left(u'w'-u''w \right)'+\mu g'(x)=0 \,.
\eeq
Integrating over $(\xi,1)$,
\[
u'(1)w'(1)-u'(\xi )w'(\xi )+\mu \int _{\xi}^1 g'(x) \, dx=0 \,.
\]
All three terms on the left are non-negative, and the second one is positive, which results in a contradiction.
\epf

\begin{lma}\lbl{lma:5}
Assume that $f(u) \in C(\bar R_+)$, and $g(x)$ satisfies the conditions (\ref{2}) and (\ref{3}). Let  $u(x)$ be a positive and singular solution of (\ref{1}), with $u'(1)<0$, and $w(x)>0$ a solution of (\ref{10}). Then
\[
\int _0^1 f(u)w \, dx>0 \,.
\]
\end{lma}

\pf
By (\ref{11}), the function $u'w'-u''w$ is non-increasing on $(0,1)$, and hence
\[
u'w'-u''w \geq u'(1)w'(1) \,.
\]
Integrating over $(0,1)$, and expressing $u''$ from (\ref{1}), gives
\[
2\int _0^1 w \left( \la f(u)-\mu g(x) \right) \, dx \geq u'(1)w'(1)>0 \,,
\]
which implies the lemma.
\epf

\begin{thm}\lbl{thm:1}
Assume that $f(u) \in C^1(\bar R_+)$, and $g(x)$ satisfies the conditions (\ref{2}) and (\ref{3}). Then
positive solutions of the problem (\ref{1}) can be continued globally  either in $\la$ or in $\mu$, on smooth solution curves, so long as $u'(1)<0$.
\end{thm}

\pf
At any non-singular solution of (\ref{1}), the implicit function theorem applies (see e.g., L. Nirenberg \cite{N}, or P. Korman  \cite{K} for more details), while at the singular solutions the Crandall-Rabinowitz \cite{CR} bifurcation theorem applies, with Lemma \ref{lma:5} verifying its crucial ``transversality condition", see e.g., P. Korman  \cite{K} (or \cite{KLO}, \cite{OS}) for more details. In either case we can always continue the solution curves.
\epf

\begin{thm}\lbl{thm:2}
{\bf (i}) Assume that $f(u) \in C^2[0,\infty)$ is concave. Then only turns to the right are possible in the $(\la ,u(0))$ plane, when solutions are continued in $\la$, and only turns to the left are possible in the $(\mu ,u(0))$ plane, when solutions are continued in $\mu$.
\medskip

\noindent
{\bf (ii}) Assume that $f(u) \in C^2[0,\infty)$ is convex. Then only turns to the left are possible in the $(\la ,u(0))$ plane, when solutions are continued in $\la$, and only turns to the right are possible in the $(\mu ,u(0))$ plane, when solutions are continued in $\mu$.
\end{thm}

\pf
Assume that when continuing in $\la$, we encounter  a critical point $(\la _0,u_0)$, i.e., the problem (\ref{10}) has a non-trivial solution $w(x)>0$. By Lemma \ref{lma:5}, the Crandall-Rabinowitz \cite{CR} bifurcation theorem applies. This theorem implies that the solution set near $(\la _0,u_0)$ is given by a curve $(\la (s),u(s))$ for $s \in (-\delta,\delta)$, with $\la (s)=\la _0+\frac12 \la '' (0)s^2+o(s^2)$, and 
\[
\la '' (0)=-\la _0 \frac{\int_{-1}^1 f''(u)w^3 \,dx}{\int_{-1}^1 f(u)w \,dx} \,,
\]
see e.g., \cite{K}.
When $f(u)$ is concave (convex), $\la '' (0)$ is positive (negative), and a turn to the right (left) occurs on the solution curve.
\medskip

If a critical point $(\mu _0,u_0)$ is encountered when continuing in $\mu $, the Crandall-Rabinowitz \cite{CR} bifurcation theorem
implies that the solution set near $(\mu _0,u_0)$ is given by $(\mu (s),u(s))$ for $s \in (-\delta,\delta)$, with $\mu (s)=\mu _0+\frac12 \mu '' (0)s^2+o(s^2)$, and 
\[
\mu '' (0)=\la \frac{\int_{-1}^1 f''(u)w^3 \,dx}{\int_{-1}^1 g(x)w \,dx} \,,
\]
see e.g., \cite{SS}.
When $f(u)$ is concave (convex), $\mu '' (0)$ is negative (positive), and a turn to the left (right) occurs on the solution curve.
\epf

\section{Numerical computation of the solution curves }
\setcounter{equation}{0}

In this section we present  computations of  the global  curves of positive solutions for the problem
(\ref{1}), which are based on  our paper \cite{K1}.
We assume that the conditions of Lemma \ref{lma:2} hold, so  that $\al \equiv u(0)$ is a global parameter. We think of the parameter $\mu$ as secondary, and we  begin with the problem (\ref{4}) (i.e., we set $\mu g(x)=h(x)$).  Since any positive solution $u(x)$ is an even function, we shall compute it on the half-interval $(0,1)$, by solving
\beq
\lbl{n2}
u''+\la f(u)-h(x)=0 \s\s \mb{for $0<x<1$,} \s\s u'(0)=u(1)=0  \,.
\eeq
A standard approach to numerical computations involves  continuation in $\la$ by using the predictor-corrector methods, see e.g., E.L. Allgower and K. Georg \cite{A}. These methods are well developed, but not easy to implement, because the solution curve $u=u(x,\la)$ may consist of several parts, each having multiple turns. Here $\la $ is a local parameter, but not a global one.
\medskip

Since  $\al = u(0)$ is a global parameter, we shall compute the solution curve $(\la,u(0))$ of (\ref{n2}) in the form $\la=\la (\al)$, with $\al=u(0)$. If we solve the initial value problem
\beq
\lbl{n3}
u'' + \la f(u)-h(x)=0, \s\s u(0)=\al,  \s\s u'(0)=0  \,,
\eeq
then we need to find $\la$, so that $u(1)=0$, in order to obtain the solution of (\ref{n2}), with $u(0)=\al$. Rewrite the equation (\ref{n3}) in the integral form
\[
u(x)=\al -\la \int_0^x (x-t) f(u(t)) \,dt+\int_0^x (x-t) h(t) \,dt \,,
\]
and then the equation for $\la$ is
\beq
\lbl{n4}
F(\la) \equiv u(1)= \al -\la \int_0^1 (1-t) f(u(t)) \,dt +\int_0^1 (1-t) h(t) \,dt =0 \,.
\eeq
We solve this equation by using Newton's method
\[
\la _{n+1}=\la _{n}-\frac{F(\la _{n})}{F \,'(\la _{n})} \,.
\]
We have
\[
F(\la _{n})=\al -\la _n  \int_0^1 (1-t) f(u(t,\la _n)) \,dt +\int_0^1 (1-t) h(t) \,dt\,,
\]
\[
F'(\la _{n})=- \int_0^1 (1-t) f(u(t,\la _n)) \,dt -\la _n  \int_0^1 (1-t) f'(u(t,\la _n)) u_{\la} \,dt\,,
\]
where $u(x,\la _n)$ and $u_{\la}$ are respectively the solutions of 
\beq
\lbl{n5}
u'' + \la _n f(u)-h(x)=0, \s u(0)=\al,  \s u'(0)=0  \,,
\eeq
\beq
\lbl{n6}
\s\s\;  u_{\la}'' + \la _n  f'(u(x,\la _n))u_{\la}+f(u(x,\la _n))=0, \s u_{\la}(0)=0,  \; u_{\la}'(0)=0  \,.
\eeq
(As we vary $\la$, we keep $u(0)=\al$ fixed, that is why $u_{\la}(0)=0$.)
This method is very easy to implement. It  requires  repeated solutions of the initial value problems (\ref{n5}) and (\ref{n6}) (using the NDSolve command in {\em Mathematica}).
\medskip

\noindent
{\bf Example} Using  {\em Mathematica} software, we have computed the solution curves in $(\la ,u(0))$ plane for the problem
\beq
\lbl{n10}
\s\s\s u''+\la u(10-2u)-\mu (1+0.2x^2)=0, \; -1<x<1, \; u(-1)=u(1)=0 \,,
\eeq
at $\mu=0.9$, $\mu=1.5$ and $\mu=2.2$. Results are presented in the Figure $1$. (The curve in the middle corresponds to $\mu=1.5$.)

\begin{figure}
\begin{center}
\scalebox{1.1}{\includegraphics{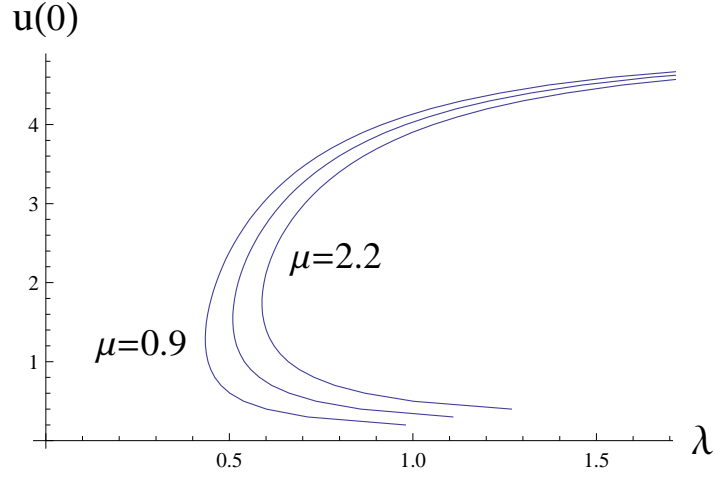}}
\end{center}
\caption{ The   curve of positive solutions for  the problem (\ref{n10}) at $\mu=0.9$, $\mu=1.5$ and $\mu=2.2$ }
\end{figure}
\medskip

We now discuss numerical continuation of solutions  in the secondary parameter. We consider again
\beq
\lbl{n7}
u''+\la f(u)-\mu g(x)=0, \s -1<x<1, \s u(-1)=u(1)=0 \,,
\eeq
and assume $\la$ to be fixed, and we continue the solutions in $\mu$.
By the Corollary \ref{cor:1}, $\al = u(0)$ is a global parameter, and we shall compute the solution curve $(\mu ,u(0))$ of (\ref{n2}) in the form $\mu=\mu (\al)$.
As before, to find $\mu=\mu (\al)$, we need to solve
\[
F(\mu) \equiv u(1)=\al -\la   \int_0^1 (1-t) f(u(t)) \,dt +\mu \int_0^1 (1-t) g(t) \,dt =0\,.
\]
We solve this equation by using Newton's method
\[
\mu _{n+1}=\mu _{n}-\frac{F(\mu _{n})}{F \,'(\mu _{n})} \,,
\]
with
\[
F(\mu _{n})=\al -\la \int_0^1 (1-t) f(u(t,\mu _n)) \,dt +\mu _n \int_0^1 (1-t) g(t) \,dt\,,
\]
\[
F'(\mu _{n})= -\la   \int_0^1 (1-t) f'(u(t,\mu _n)) u_{\mu} \,dt+ \int_0^1 (1-t) g(t) \,dt\,,
\]
where $u(x,\mu _n)$ and $u_{\mu}$ are respectively the solutions of 
\[
u'' + \la f(u)-\mu _n g(x)=0, \s u(0)=\al,  \s u'(0)=0  \,,
\]
\[
\s\s\;  u_{\mu}'' + \la   f'(u(x,\mu _n))u_{\mu}-g(x)=0, \s u_{\mu}(0)=0,  \; u_{\mu}'(0)=0  \,.
\]
(As we vary $\mu$, we keep $u(0)=\al$ fixed, that is why $u_{\mu}(0)=0$.)

\begin{figure}
\begin{center}
\scalebox{.7}{\includegraphics{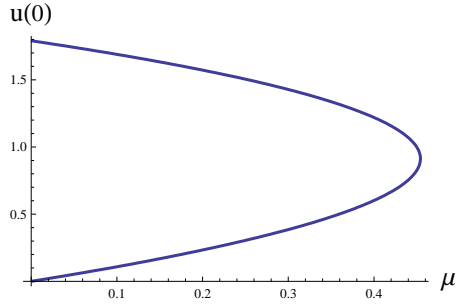}}
\end{center}
\caption{ The  curve of positive solutions for  the problem (\ref{22}) }
\end{figure}
\medskip

\noindent
{\bf Example}
We have continued in $\mu$ the positive solutions of
\beq
\lbl{22}
\s\s u''+ u(4-u)-\mu (1+x^2)=0, \s -1<x<1, \s u(-1)=u(1)=0 \,.
\eeq
The curve of positive solutions is given in Figure $2$. 
\medskip

\noindent
{\bf Example}
We have continued in $\mu$ the positive solutions of
\beq
\lbl{20}
\s\s u''+2.4 u(4-u)-\mu (1+x^2)=0, \s -1<x<1, \s u(-1)=u(1)=0 \,.
\eeq
The curve of positive solutions is given in Figure $3$. At $\mu \approx 2.28634$, the solutions become sign changing, negative near $x=\pm 1$.

\begin{figure}
\begin{center}
\scalebox{.7}{\includegraphics{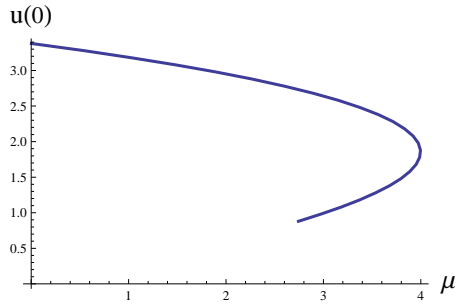}}
\end{center}
\caption{ The  curve of positive solutions for  the problem (\ref{20}) }
\end{figure}

\section{Diffusive logistic equation with  harvesting }
\setcounter{equation}{0}
Recall that the eigenvalues of 
\[
u''+\la u=0, \s -1<x<1, \s u(-1)=u(1)=0
\]
are $\la _n=\frac{n^2 \pi^2}{4}$, and in particular $\la _1=\frac{ \pi^2}{4}$, $\la _2=\pi ^2$.
\medskip

We consider positive solutions of 
\beq
\lbl{30}
 u''+\la u(1-u)-\mu =0, \s -1<x<1, \s u(-1)=u(1)=0 \,,
\eeq
with positive parameters $\la$ and $\mu$.
It is easy to see that no positive solutions exist if $\la \leq \la _1=\frac{ \pi^2}{4}$, and by maximum principle any positive solution satisfies $0<u(x)<1$. The following result gives a complete description of the set of positive solutions. 

\begin{thm}\lbl{thm:3}
For any fixed $\mu$ the set of positive solutions of (\ref{30}) is a parabola-like curve opening to the right in $(\la ,u(0))$ plane (the $\la$-curves). The upper branch continues for all $\la$ after the turn, while the solutions on the lower branch become sign-changing after some $\la=\bar \la$ ($u_x(\pm 1,\bar \la )=0$, see Figure $1$). For any fixed $\la$ the set of positive solutions of (\ref{30})  is a parabola-like curve opening to the left in $(\mu ,u(0))$ plane (the $\mu$-curves). Different $\la$-curves (and different $\mu$-curves) do not intersect. The $\la$-curves and the $\mu$-curves share the turning points. Namely, if at $\mu =\mu _0$, the $\la$-curve turns at the point $(\la _0,\al)$, then at $\la =\la _0$, the $\mu$-curve turns at the point $(\mu _0,\al)$.
\medskip

If $\la \in (\la _1, \la _2]$, then the $\mu$-curve joins the point $(0,\mu _1)$, with some $\mu _1>0$,  to $(0,0)$, with exactly one turn to the left at some $\mu _0$ (as in Figure $2$). If $\la > \la _2$, then the $\mu$-curve joins the point $(0,\mu _1)$, with some $\mu _1>0$,  to some point $(\bar \mu>0,\al >0)$, with exactly one turn to the left at some $\mu _0>\bar \mu$ (as in Figure $3$). Solutions on the lower branch become sigh-changing for $\mu<\bar \mu$.
\end{thm}

We remark that in case $\la \in (\la _1, \la _2]$, a more general result was given in J. Shi \cite{S}, by a more involved method.

The proof will depend on several lemmas, which we state for a more general problem 
\beq
\lbl{31}
 u''+\la f(u)-\mu =0, \s -1<x<1, \s u(-1)=u(1)=0 \,.
\eeq

\begin{lma}\lbl{lma:30}
Assume that $f(u) \in C^1(\bar R_+)$ satisfies $f(0)=0$, and $f(u(x))>0$ for any positive solution of (\ref{31}), for all $x \in (-1,1)$. Assume that $u(x,\la)$ arrives at the point $\la _0$ where the positivity of solutions is lost (i.e., $u_x(\pm 1,\la _0)=0$) with the maximum value $u(0,\la)$ decreasing along the solution curve. Then the positivity is lost forward in $\la$ at $\la _0$. (I.e., $u(x,\la)>0$ for all $x \in (-1,1)$ if $\la<\la _0$, and $u(x,\la)$ is sign-changing for $\la>\la _0$.)
\end{lma}

\pf
Since for positive solutions $u_x(x,\la _0)<0$ for $x \in (0,1)$, the only way for solutions to become sign-changing is to have $u_x(\pm 1,\la _0)=0$. Assume, on the contrary, that positivity is lost backward in $\la$. Then by our assumption
\beq
\lbl{32}
u_{\la}(0,\la _0) \geq 0 \,.
\eeq
Differentiating the equation (\ref{31}) in $\la$, we have
\beq
\lbl{33}
u_{\la }''+\la f'(u)u_{\la } =-f(u), \s -1<x<1, \s u_{\la }(-1)=u_{\la }(1)=0 \,.
\eeq
Differentiating the equation (\ref{31}) in $x$, gives
\beq
\lbl{34}
u_{x }''+\la f'(u)u_{x } =0 \,.
\eeq
Combining the equations (\ref{33}) and (\ref{34}),
\[
\left(u_{\la }'u'-u_{\la } u'' \right)'=-f(u)u_x>0 \,, \s \mbox{for $x \in (0,1)$} \,.
\]
It follows that the function $q(x) \equiv u_{\la }'u'-u_{\la } u''$ is increasing, with $q(0) \geq 0$ by (\ref{32}), and $q(1)=0$, a contradiction.
\epf

\begin{lma}\lbl{lma:31}
Assume that $f(u) \in C^1(\bar R_+)$ satisfies $f(0)=0$, and $f(u(x))>0$ for any positive solution of (\ref{31}), for all $x \in (-1,1)$. Assume that $u(x,\mu)$ arrives at the point $\mu _0$ where the positivity is lost (i.e., $u_x(\pm 1,\mu _0)=0$) with the maximum value $u(0,\mu)$ decreasing along the solution curve. Then the positivity is lost backward in $\mu$ at $\mu _0$. (I.e., $u(x,\mu)>0$ for all $x \in (-1,1)$ if $\mu > \mu _0$, and $u(x,\mu)$ is sign-changing for $\mu<\mu _0$.)
\end{lma}

\pf
Assume, on the contrary, that positivity is lost forward in $\mu$. Then by our assumption
\beq
\lbl{36}
u_{\mu}(0,\mu _0) \leq 0 \,.
\eeq
Differentiating the equation (\ref{31}) in $\mu$, we have
\beq
\lbl{37}
u_{\mu }''+\la f'(u)u_{\mu } =1, \s -1<x<1, \s u_{\mu }(-1)=u_{\mu }(1)=0 \,.
\eeq
Combining the equations (\ref{34}) and (\ref{37}),
\[
\left(u_{\mu }'u'-u_{\mu } u'' \right)'=u_x<0 \,, \s \mbox{for $x \in (0,1)$} \,.
\]
It follows that the function $r(x) \equiv u_{\mu }'u'-u_{\mu } u''$ is decreasing, with $r(0) \leq 0$ by (\ref{36}), and $r(1)=0$, a contradiction.
\epf

\begin{lma}\lbl{lma:33}
Assume that $f(u) \in C^1[0,\infty)$, $f(0)=0$, and $f'(u)$ is a decreasing function for $u>0$. Then for any $\la >0$ there is at most one solution pair $(\mu,u(x))$, with $u'(\pm 1)=0$.
\end{lma}

\pf
Assume, on the contrary, that there are two solution pairs $(\mu _1,u(x))$ and $(\mu _2,v(x))$,with $\mu _2>\mu _1$, satisfying 
\beq
\lbl{38}
 u''+\la f(u)-\mu _1 =0, \s -1<x<1, \s
 u(\pm 1)=u'(\pm 1)= 0 \,,
\eeq
\beq
\lbl{39}
 v''+\la f(v)-\mu _2 =0, \s -1<x<1, \s
 v(\pm 1)=v'(\pm 1)= 0 \,.
\eeq
Since $v''(1)=\mu _2>\mu _1=u''(1)$, we have $v(x)>u(x)$ for $x$ close to $1$. Two cases are possible.
\medskip

\noindent
{\bf (i)} $v(x)>u(x)$ for $x \in (0,1)$.  Differentiating the equations (\ref{38}) and (\ref{39}), we get
\[
u_x''+\la f'(u)u_x=0, \s \mbox{$u_x<0$ on $(0,1)$}, \s u_x(0)=u_x(1)=0 \,,
\]
\[
v_x''+\la f'(v)v_x=0, \s \mbox{$v_x<0$ on $(0,1)$}, \s v_x(0)=v_x(1)=0 \,.
\]
Since $f'(u)>f'(v)$, we have a contradiction by Sturm's comparison theorem.
\medskip

\noindent
{\bf (ii)} There is  $\xi \in (0,1)$ such that $v(x)>u(x)$ for $x \in (\xi,1)$, while $v(\xi)=u(\xi)$ and $u'(\xi) \leq v'(\xi)<0$. Multiply the equation (\ref{38}) by $u'$, and integrate over $(\xi,1)$ (with $F(u)=\int_0^u f(t) \, dt$)
\[
-\frac12 {u'}^2(\xi)-\la F(u(\xi))+\mu _1 u(\xi)=0 \,.
\]
Similarly, from (\ref{39})
\[
-\frac12 {v'}^2(\xi)-\la F(u(\xi))+\mu _2 u(\xi)=0 \,.
\]
Subtracting
\[
(\mu _2 -\mu _1)u(\xi)=\frac12 \left({v'}^2(\xi)-{u'}^2(\xi) \right) \,.
\]
The quantity on the left is positive, while the one on the right is non-positive, a contradiction.
\epf

\begin{lma}\lbl{lma:34}
For any $\al \in (0,\frac34 )$ there exists a unique pair $(\bar \la,\bar  \mu)$, with $\bar \la >\la _2$ and $\bar  \mu>0$,  and a positive solution of   (\ref{30}) with $u(0)=\al$ and $u'(\pm 1)=0$.
Moreover, if $\bar \mu \ra 0$, then $ \bar \la \downarrow \la_2=\pi^2$.
\end{lma}

\pf
Multiplying the equation (\ref{30}) by $u'$, we see that the solution with $u'(\pm 1)=0$ satisfies
\beq
\lbl{40}
\frac12 {u'}^2+\bar \la \left( \frac12 u^2-\frac13 u^3 \right)-\bar \mu u=0 \,.
\eeq
Evaluating this at $x=0$
\beq
\lbl{41}
\bar \la \left( \frac12 \al -\frac13 \al ^2 \right)=\bar \mu \,.
\eeq
We also express from (\ref{40})
\[
\frac{du}{dx}=- \sqrt{ 2 \bar \mu u-\bar \la \left(  u^2-\frac23 u^3 \right)} \,, \s  \mbox{for $x \in (0,1)$} \,.
\]
We express $\bar \mu$ from (\ref{41}), separate the variables and integrate, getting
\[
\int_0^{\al } \frac{du}{\sqrt{ \left(\al-\frac23 \al ^2 \right)  u- \left(  u^2-\frac23 u^3 \right)}}=\sqrt{\bar \la} \,.
\]
Setting here $u=\al v$, we express
\beq
\lbl{42}
\bar \la=\left( \int_0^{1 } \frac{dv}{\sqrt{ \left(1-\frac23 \al  \right)  v- \left(  v^2-\frac23 \al v^3 \right)}} \right)^2 \,.
\eeq
The formulas (\ref{42}) and (\ref{41}) let us compute $\bar \la$, and then $\bar \mu$, for any $\al \in (0,\frac34 )$. (For  $\al \in (0,\frac34 )$, the quantity inside the square root in (\ref{42}), which is $v(1-v) \left[1  -\frac23 \al (1+v) \right]$, is positive for all $v \in (0,1)$.)
\medskip

If $\bar \mu \ra 0$, then from (\ref{41}), $\al \ra 0$ (recall that $\la >\la _1$), and then from (\ref{42})
\[
\bar \la \downarrow \left( \int_0^{1 } \frac{dv}{\sqrt{   v-   v^2}} \right)^2 =\pi ^2\,, 
\]
completing the proof. 
\epf

\noindent
{\bf Proof of the Theorem \ref{thm:3}}
It is easier to understand the $\la$-curves, so we assume first that $\mu$ is fixed. It is well known that for $\la$ large enough the problem (\ref{30}) has a positive stable (``large") solution, with $u(0,\la)$ increasing in $\la$ (see e.g., \cite{SS}). Let us continue this solution for decreasing $\la$. This curve does not  continue to $\la$'s $\leq \la _1$, and it cannot become sign-changing while continued to the left, by Lemma \ref{lma:30}, hence a turn to the right must occur. After the turn, standard arguments imply that solutions develop zero slope at $\pm 1$, and become sign-changing for $\la >\bar \la _{\mu}$, see e.g., \cite{K}. By Theorem \ref{thm:2}, exactly one turn occurs on each $\la$-curve, and by Lemma \ref{lma:34}, $\inf _{\mu} \bar \la _{\mu}=\la _2=\pi ^2$.
\medskip

Turning to the $\mu$-curves, for any fixed $\tilde \la >\la _1$ we can find a positive solution on the curve $\mu=0$ (the curve that bifurcates from the trivial solution at $\la =\la _1$). We now slide down from this point in the $(\la, u(0))$ plane, by varying $\mu$. As we increase $\mu$ (keeping $\tilde \la$ fixed), we slide to different $\la$-curves. At some $\mu$ we reach a $\la$-curve which has its turn at $\la =\tilde \la$. After that point, $\mu$ begins to decrease on the $\la$-curves. If $\tilde \la \in (\la _1,\la _2]$, we slide all the way to $\mu=0$, by Lemma \ref{lma:34}. Hence, the $\mu$-curve at $\tilde \la$ is as in Figure $2$. In case $\tilde \la >\la _2$, by Lemma \ref{lma:34}, we do not slide all the way to  $\mu=0$, and hence the $\mu$-curve at $\tilde \la$ is as in Figure $3$. By Lemma \ref{lma:33}, this curve exhausts the set of all positive solutions of (\ref{30}) (any other solution would lie on a solution curve with no place to go, when continued in $\mu$).
\epf

\noindent
{\bf Remark} Our results also imply that the $\mu$-curves described in the Theorem \ref{thm:3} continue without turns for all $\mu <0$. 
(Observe that Lemma \ref{lma:4} holds for autonomous problems, regardless of the sign of $\mu$, see e.g., \cite{K}.) Negative $\mu$'s correspond to ``stocking" of fish, instead of ``fishing". In Figure $4$ we present the solution curve of the problem
\beq
\lbl{46}
 u''+6u(1-u)-\mu =0, \s -1<x<1, \s u(-1)=u(1)=0 \,.
\eeq
Observe that $u(0)>1$ for  $\mu<\mu _0$, for some  $\mu _0<0$.
\medskip

For the non-autonomous version of the fishing problem
\beq
\lbl{50}
 u''+\la u(1-u)-\mu g(x)=0, \s -1<x<1, \s u(-1)=u(1)=0 
\eeq
we were not able to extend any of the above lemmas. Still it appears easier to understand the $\la$-curves first. Here we cannot rule out the possibility of the $\la$-curves losing their positivity backward, and thus never making a turn to the right. 
\medskip

We prove next that a turn to the right does occur for solution curves of (\ref{50}) that are close to the curve bifurcating from zero at $\la =\la _1$. Let $\bar \la _{\mu}$ be the value of $\la$, at which positivity is lost for a given $\mu$ (so that $u_x(\pm 1,\bar \la _{\mu})=0$). We claim that $\inf _{\mu>0} \bar \la _{\mu}>\la _1$. Indeed, assuming otherwise, we can find a sequence $\{\mu _n \} \ra 0$, with $\bar \la _{\mu _n} \ra \la _1$. By a standard argument, $\ds \frac{u(x,\bar \la _{\mu_n})}{u(0,\bar \la _{\mu_n})} \ra w(x)>0$, where 
\[
w''+\la _1w=0, \s w(\pm 1)=w'(\pm1)=0 \,,
\]
which is not possible. It follows that for a fixed $\la \in (\la _1, \inf _{\mu>0} \bar \la _{\mu}]$, the $\mu$-curve for (\ref{50}) is as in Figure $2$. (Notice that this also implies that the $\la$-curves for small $\mu$ do turn.) A similar result for general PDE's was proved in S.  Oruganti, J. Shi, and R. Shivaji \cite{SS}. In case $\la> \inf _{\mu>0} \bar \la _{\mu}$, the $\mu$-curves are different, although we cannot prove in general that they are as in Figure $3$. (For example, we cannot rule out the possibility that the $\mu$-curves consist of several pieces.)
\begin{figure}
\begin{center}
\scalebox{.7}{\includegraphics{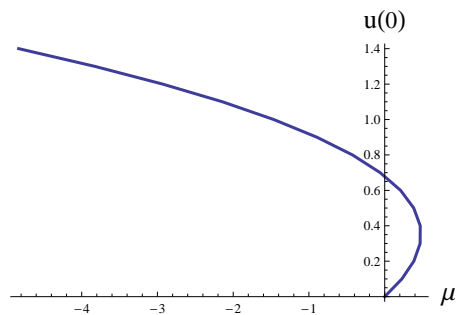}}
\end{center}
\caption{ The  solution curve for  the problem (\ref{46}) }
\end{figure}


\begin{thebibliography}{99}

\bibitem{A}
E.L. Allgower and K. Georg, Numerical Continuation Methods. An Introduction. Springer Series in Computational Mathematics, {\bf 13}. Springer-Verlag, Berlin, (1990). 
\vspace{-0.2cm}



\bibitem{CR} 
M.G. Crandall and P.H. Rabinowitz,  Bifurcation, perturbation of simple
eigenvalues and linearized stability, {\em Arch. Rational Mech. Anal.}, {\bf 52}, 161-180 (1973).
\vspace{-0.2cm}

\bibitem{C}
D.G. Costa, P. Dr\'{a}bek, and H. Tehrani, Positive solutions to semilinear elliptic equations with logistic type nonlinearities and constant yield harvesting in $R^n$, {\em Comm. Partial Differential Equations} {\bf 33}, 1597-1610  (2008). 
\vspace{-0.2cm}

\bibitem{GNN} 
B. Gidas, W.-M. Ni and L. Nirenberg, Symmetry and related properties via
the maximum principle, {\em Commun. Math. Phys.} {\bf 68}, 209-243 (1979).
\vspace{-0.2cm}


\bibitem{G}
P. Gir\~{a}o, and H. Tehrani,  Positive solutions to logistic type equations with harvesting, {\em  J. Differential Equations} {\bf 247}, no. 2, 574-595  (2009).
\vspace{-0.2cm}

\bibitem{H}
S.P. Hastings and J.B.  McLeod,  Classical Methods in Ordinary Differential Equations. With applications to boundary value problems. Graduate Studies in Mathematics, {\bf 129}. American Mathematical Society, Providence, RI (2012).
\vspace{-0.2cm}

\bibitem{W1} 
K.C. Hung and S.H. Wang, Classification and evolution of bifurcation curves for a multiparameter p-Laplacian Dirichlet problem, {\em Nonlinear Anal.} {\bf  74}, no. 11, 3589-3598 (2011).
\vspace{-0.2cm}
 
\bibitem{W2} 
K.C. Hung and S.H. Wang, Bifurcation diagrams of a p-Laplacian Dirichlet problem with Allee effect and an application to a diffusive logistic equation with predation, {\em  J. Math. Anal. Appl.} {\bf 375}, no. 1, 294-309  (2011).
\vspace{-0.2cm}

\bibitem{K14}
P. Korman, Symmetry of positive solutions for elliptic problems in one dimension, {\em Appl. Anal.} {\bf 58}, no. 3-4, 351-365  (1995).
\vspace{-0.2cm}

\bibitem{K}
P. Korman, Global Solution Curves for Semilinear Elliptic Equations, World Scientific, Hackensack, NJ (2012).
 \vspace{-0.2cm}

\bibitem{K1} 
P. Korman,  Exact multiplicity and numerical computation of solutions for two classes of non-autonomous problems with concave-convex nonlinearities, {\em Nonlinear Anal.} {\bf 93}, 226-235  (2013).
 \vspace{-0.2cm}

\bibitem{KLO}
P. Korman, Y. Li and T. Ouyang, Exact multiplicity results for boundary-value
problems with nonlinearities generalising cubic,  {\em Proc.
Royal Soc. Edinburgh, Ser. A} {\bf 126A}, 599-616 (1996).
\vspace{-0.2cm}

\bibitem{N}
L. Nirenberg, Topics in Nonlinear Functional Analysis, Courant Institute Lecture Notes, {\em Amer. Math. Soc.} (1974).
\vspace{-0.2cm}

\bibitem{SS}
S.  Oruganti, J. Shi, and R. Shivaji,  Diffusive logistic equation with constant yield harvesting. I. Steady states, {\em Trans. Amer. Math. Soc.} {\bf  354}, no. 9, 3601-3619  (2002).
\vspace{-0.2cm}

\bibitem{OS}
T. Ouyang and J. Shi,
Exact multiplicity of positive solutions for a class of semilinear problems, II, {\em J. Differential Equations} {\bf 158},  no. 1, 94-151 (1999).
\vspace{-0.2cm}

\bibitem{S}
J.  Shi, A radially symmetric anti-maximum principle and applications to fishery management models, {\em  Electron. J. Differential Equations}, No. 27, 13 pp. (electronic) (2004).
\vspace{-0.2cm}


\end{thebibliography}
\end{document}